\documentclass[12pt]{article} 
\usepackage[reqno]{amsmath} 
\usepackage{amssymb, amsthm}
\usepackage{enumerate}  
\usepackage{hhline}
\usepackage{url} 
\usepackage{graphicx}
\usepackage{float}
\usepackage{nicefrac}

\newtheoremstyle{plainsl}%
	{\topsep}
	{\topsep}
	{\slshape} 
	{}
	{\normalfont\bfseries}
	{.}
	{ }
	{}

\swapnumbers

\theoremstyle{plainsl}
\newtheorem{theorem}{Theorem}[section]
\newtheorem{lemma}[theorem]{Lemma}

\newtheorem{corollary}[theorem]{Corollary}
\newtheorem*{theorem*}{Theorem}
\newtheorem*{question}{Question}

\renewcommand\proof{\noindent\textsl{Proof. }}
\newcommand\sqr[2]{{\vbox{\hrule height.#2pt
    \hbox{\vrule width.#2pt height#1pt \kern#1pt
        \vrule width.#2pt}\hrule height.#2pt}}}

\renewcommand\qed{%
	\ifmmode\eqno\sqr53
	\else\nolinebreak\ \hfill\sqr53\medbreak\fi}

\numberwithin{equation}{section}

\newcommand\mxm{\widehat{M}}

\newcommand\Zv{{\mathbf v}}

\newcommand\opk[1]{\mathop{\mathrm{#1}}\nolimits}
\newcommand\comp[1]{{\mkern2mu\overline{\mkern-2mu#1}}}
\newcommand\diff{\mathbin{\mkern-1.5mu\setminus\mkern-1.5mu}}

\newcommand\pmat[1]{\begin{pmatrix} #1 \end{pmatrix}}

\newcommand\ones{{\bf1}}
\newcommand\rk{\opk{rk}}
\newcommand\tr{\opk{tr}}

\title{Average mixing matrix of trees}
\author{Chris Godsil \thanks{University of Waterloo, Waterloo, Canada. \protect\url{{cgodsil,john.sinkovic}@uwaterloo.ca}}, Krystal Guo\thanks{Universit\'{e} libre de Bruxelles, Brussels, Belgium. \protect\url{guo.krystal@gmail.com}} \thanks{This work was done when K. Guo was a postdoctoral fellow at the University of Waterloo.}, John Sinkovic\footnotemark[1]}

\begin{document}

\maketitle

\begin{abstract}
We investigate the rank of the average mixing matrix of trees, with all eigenvalues distinct. The rank of the average mixing matrix of a tree on $n$ vertices with $n$ distinct eigenvalues is upper-bounded by $\nicefrac{n}{2}$. Computations on trees up to $20$ vertices suggest that the rank attains this upper bound most of the times. We give an infinite family of trees whose average mixing matrices have ranks which are bounded away from this upper bound. We also give a lower bound on the rank of the average mixing matrix of a tree. 
 \end{abstract}

\section{Introduction}
\label{sec:intro}
We investigate the rank of the average mixing matrix of continuous time quantum walks on trees with simple eigenvalues.

The continuous-time quantum walk is an universal computational primitive \cite{ChildsUniversalQComputation}. Introduced in \cite{FarhiGutmann}, many properties of quantum walks have been studied. Some topics include state transfer \cite{KayReviewPST,GodsilStateTransfer12,CoutinhoGodsilGuoVanhove2,VinetZhedanovHowTo} and uniform mixing \cite{AdaChanComplexHadamardIUMPST,GodsilMullinRoy,TamonAdamczakUniformMixingCycles}. 

Let $X$ be a graph and let $A$ be the adjacency matrix of $X$. The \textsl{transition matrix} of a continuous-time quantum walk on $X$ is a matrix-valued function in time, denoted $U(t)$, and is given as follows:
\[ U(t) = \exp(itA)
\]
Unlike a classical random walk on a connected graph, a continuous quantum walk does not reach a stationary distribution. For example, the transition of the complete graph $K_n$ is a periodic function of $t$. The average mixing matrix is, intuitively, a distribution that the quantum walk adheres to, on average, over time, and thus may be thought of as a replacement for a stationary distribution. The average mixing matrix has been studied in \cite{GodsilAverageMixing,TamonAdamzackAverageMixing, CoGdGuZh17}.  We will now proceed with a few preliminary definitions. 

The mixing matrix $M(t)$ of the quantum walk, given by
\[
	M(t) = U(t) \circ \comp{U(t)}.
\]
where $\circ$ denotes the Schur (also known as the Hadamard or element-wise product) of two matrices. The \textsl{average mixing matrix}, denoted $\mxm$, is defined as follows: 
\begin{equation}
	\label{eq:mxmdef}
	\mxm = \lim_{T\to\infty}\frac1T \int_0^T M(t)\,dt.
\end{equation}

A tree is a connected graph with no cycles. The class of trees is a natural class of graphs to study. The rank of the average mixing matrix of a tree on $n$ vertices with $n$ distinct eigenvalues is at most $\nicefrac{n}{2}$, see \cite{CoGdGuZh17}. We computed the rank of the average mixing matrix of all trees up to $20$ vertices; the computational results suggest that the rank of $\mxm$ of a tree with all simple eigenvalues attains the upper bound most of the time. We give an infinite family of trees whose average mixing matrices have ranks which are bounded away from this upper bound.

\begin{theorem*}
For every positive real number $c$, there exists a tree $T$ on $n$ vertices with simple eigenvalues such that 
\[
\left\lceil \frac{n}{2} \right\rceil -  \rk(\mxm(T)) > c.
\]
\end{theorem*}

We then give the following lower bound on the rank of the average mixing matrix of a tree:

\begin{theorem*}
If $T$ is a tree with simple eigenvalues on at least four vertices and is not isomorphic to $P_4$, then $\rk(\mxm(T))\geq 3$.
\end{theorem*}

We give some preliminary results about the average mixing matrix in Section \ref{sec:prelim}. We also find the ranks and traces of the average mixing matrices of star graphs, as an example. In Section \ref{sec:rootprod}, we define the rooted product of graphs and we find the average mixing matrix of a graph $X$ with $K_2$ in terms of the eigenprojectors of $X$, which we use in Section \ref{sec:upperbd} to give a family of trees whose average mixing matrices have ranks which are upper-bounded away from the maximum. In Section \ref{sec:lowerbd}, we prove the lower bound on the ranks of the average mixing matrices of trees which are not isomorphic to the path on four vertices. 

\section{Preliminaries and examples}\label{sec:prelim}

Let $X$ be a graph on $n$ vertices and let $A$ be the adjacency matrix of $X$. Let $\theta_1, \ldots, \theta_d$ be the distinct eigenvalues of $A$ and, for $r=1,\ldots,d$, let $E_r$ be the idempotent projection onto the $\theta_r$ eigenspace of $A$; the spectral decomposition of $A$ is as follows:
\[
A = \sum_{r=1}^d \theta_r E_r. 
\]
The following is an important theorem as it allows us to understand the average mixing matrix. 

\begin{theorem}\cite{GodsilAverageMixing}
Let $X$ be a graph and let $A$ be the adjacency matrix of $X$. Let $A = \sum_{r=1}^d \theta_r E_r$ be the spectral decomposition of $A$. The mixing matrix of $X$ is \[
\mxm = \sum_{r=1}^d E_r \circ E_r . \qed \] 
\end{theorem}

As an example, we find the average mixing matrix of the star and deduce that it always has full rank. 

\begin{lemma}
The average mixing matrix of the star graph $K_{1,n}$ has full rank and has 
\[
\tr(\mxm) = \begin{cases}\frac{2n^2 + \frac{5}{2}n}{(2n-1)^2} , & \text{if }n \text{ is even;} \\
\frac{2n^2 + 3n}{(2n-1)^2} , & \text{if }n \text{ is odd.} \end{cases}
\]
\end{lemma}

\proof We will proceed by finding an explicit expression for 
$\mxm(K_{1,n})$. Observe that the eigenvalues of $K_{1,n}$ are $\pm \sqrt{n}$ and $0$ with multiplicity $n-1$. We give the spectral idempotent matrices  as follows:
\[
E_{\sqrt{n}} = \frac{n}{2n-1} \pmat{1 & \frac{1}{\sqrt{n}} \ones^T \\
\frac{1}{\sqrt{n}} \ones & \frac{1}{n}J} , \quad
E_{-\sqrt{n}} = \frac{n}{2n-1} \pmat{1 & -\frac{1}{\sqrt{n}} \ones^T \\
-\frac{1}{\sqrt{n}} \ones & \frac{1}{n}J}
\]
and
\[
E_0 = \frac{1}{2} \pmat{0 & 0 \\ 0 & I -R}
\]
where $I$ is the $n\times n$ identity matrix and $R$ is the $n\times n$ matrix with all ones on the back diagonal. Note that $I-R$ always has rank $n-1$. Thus we have that
\[
\mxm(K_{1,n}) = E_{\sqrt{n}}^{\circ 2} + E_{-\sqrt{n}}^{\circ 2} + E_0^{\circ 2} = \frac{2}{(2n-1)^2}\pmat{n^2 & n\ones^T\\ n\ones & J + \frac{1}{4}I + \frac{1}{4}R}.
\]
Observe that this matrix has full rank.\qed

\begin{figure}[ht]
    \centering
    \includegraphics[scale=0.8]{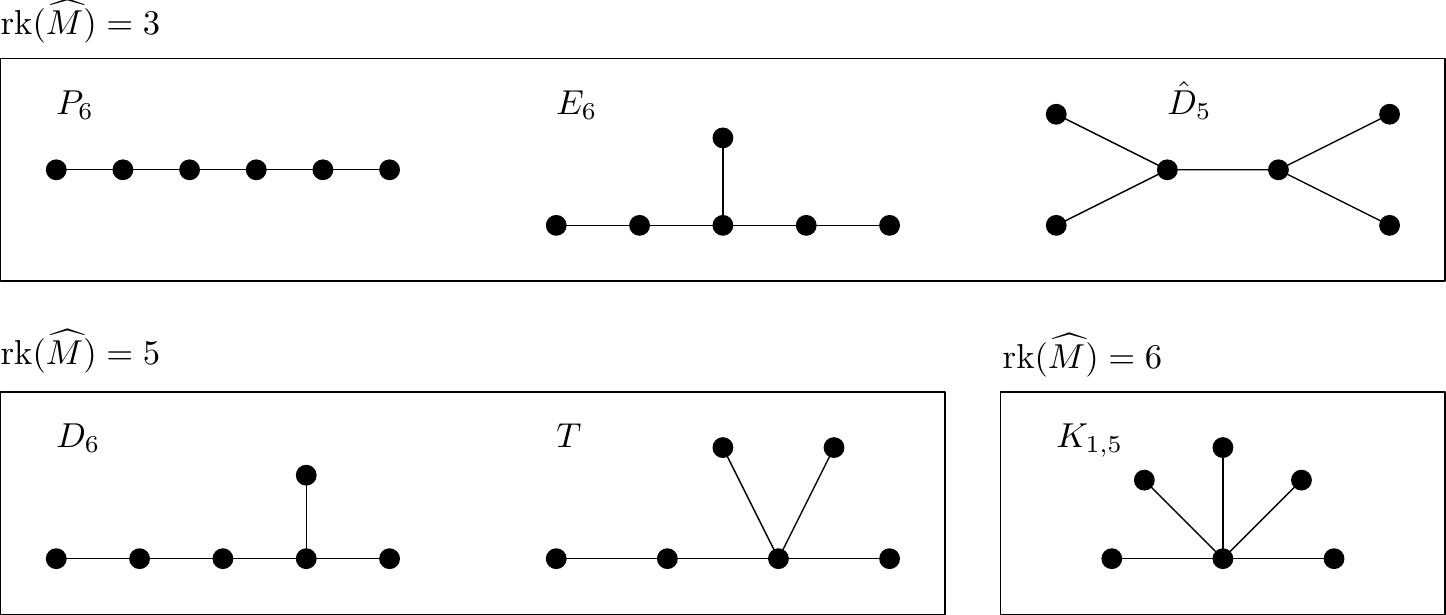}
    \caption{All trees on $6$ vertices, organized by the rank of their average mixing matrices.     \label{fig:trees6}}
\end{figure}

The path graphs $P_2$ and $P_3$ are the only trees on two and three vertices, respectively. There are six trees on six vertices; three of them have simple eigenvalues and rank $3$, two have rank $4$ and one graph (the star graph) of full rank. They are shown in Figure \ref{fig:trees6}. 

\section{Rooted products} \label{sec:rootprod}

Let $X$ be a graph with vertices $\{v_1,\ldots, v_n\}$ and let $Y$ be a disjoint union of rooted graphs $Y_1, \ldots, Y_n$, rooted at $y_1,\ldots, y_n$, respectively. The \textsl{rooted product} of $X$ and $Y$, denoted $X(Y)$, is the graph obtain by identifying $v_i$ with the root vertex of $Y_i$. The rooted product was first introduced by Godsil and McKay in \cite{GoMc78}. Let $B(X,Y)$ be the matrix given as follows:
\[
B(X,Y)_{i,j} = \begin{cases} \phi(Y_i, t), &\text{ if }i=j;\\
- \phi(Y_i \setminus y_i,t), &\text{ if }i\sim j; \\ \text{ and,}
0, & \text{ otherwise.}\end{cases} 
\]
\begin{theorem}\cite{GoMc78}
$\phi(X(Y), t) = \det(B)$. \qed 
\end{theorem}

In particular, we will consider the special case where $Y$ is a sequence of $n$ copies of $K_2$. In this case, we will write $X(K_2)$ to denote the rooted product and  the following is found in \cite{GoMc78}  as a consequence of the above theorem
\begin{equation} \label{eq:rootcharpoly}
\phi(X(K_2),t) = t^n \phi\left(H, t - \frac{1}{t}\right)
\end{equation}

The following lemma follows directly from \eqref{eq:rootcharpoly}. 
\begin{lemma} \label{lem:rootedsimp}
Let $X$ be a graph with simple eigenvalues. The rooted product $X(K_2)$ also has simple eigenvalues which are roots of 
\[
t^2 - \lambda t - 1 =0,
\]
for $\lambda$ an eigenvalue of $X$. \qed
\end{lemma} 

\begin{lemma}
Let $X$ be a graph with simple eigenvalues and let $\{\Zv_1,\ldots,\Zv_n\}$ be an orthonormal eigenbasis of $A(X)$ with eigenvalues $\lambda_1,\ldots, \lambda_n$, respectively. For $i = 1,\ldots, n$, let $\mu_i$ and $\nu_{i}$ be the two roots of $t^2 - \lambda_i t - 1 =0$. Then 
\[
\left\{\frac{1}{\sqrt{\mu_i^2 + 1}} \pmat{\mu_i \Zv_i \\  \Zv_i}, \frac{1}{\sqrt{\nu_i^2 + 1}} \pmat{\nu_i\Zv_i \\  \Zv_i} \mid i = 1,\ldots,n \right\}
\]
is an orthonormal eigenbasis for $X(K_2)$. 
\end{lemma}

\proof We may write the adjacency matrix of $X(K_2)$ as follows:
\[
A(X(K_2)) = \pmat{A(X) & I \\ I & 0} .
\]
For $i \in [n]$ and $\mu \in \{\mu_i,\nu_i\}$, we see that $\mu \neq 0$ and $\mu^2 = \lambda\mu + 1$. We obtain that 
\[ \pmat{A(X) & I \\ I & 0} \pmat{ \mu\Zv_i \\ \Zv_i} = \pmat{\lambda_i\mu\Zv_i + \Zv_i \\ \mu \Zv_i} = \pmat{\mu^2 \Zv_i \\ \mu \Zv_i} =  \mu \pmat{ \mu\Zv_i \\ \Zv_i} ,
\]
and the lemma follows. \qed

\begin{theorem} \label{thm:rootedmhat}
Let $X$ be a graph with simple eigenvalues and let $F_1,\ldots,F_n$ the orthogonal projections onto the eigenspaces of $A(X)$ with corresponding eigenvalues $\lambda_1,\ldots, \lambda_n$. Then
\[
\mxm(X(K_2)) = \pmat{\mxm(X) - N & N \\ N & \mxm(X) - N}
\]
where
\[ N =  \sum_{i=1}^n  \left(\frac{2}{\lambda_i^2 + 4}\right)(F_i\circ F_i).
\]
\end{theorem}

\proof Following the notation of the previous lemma, we see that $F_i = \Zv_i \Zv_i^T$. For $i\in [n]$, let $\mu \in \{\mu_i,\nu_i\}$. The orthogonal projection $E_{\mu} $ onto the $\mu$ eigenspace of $A(X(K_2)$ is given by 
\[
\frac{1}{\mu^2\! + \!1} \pmat{\mu \Zv_i \\  \Zv_i} \pmat{\mu \Zv_i^T &  \Zv_i^T} = \frac{1}{\mu^2\! +\! 1} \pmat{\mu^2 \Zv_i\Zv_i^T & \mu \Zv_i\Zv_i^T \\ \mu \Zv_i\Zv_i^T & \Zv_i\Zv_i^T  } = \frac{1}{\mu^2\! +\! 1} \pmat{\mu^2 F_i & \mu F_i\\ \mu F_i & F_i  }.
\]
Thus
\[
E_{\mu_i} \circ E_{\mu_i} + E_{\nu_i} \circ E_{\nu_i}\! =  \!
\pmat{\left(\!\frac{\mu_i^4 }{\left(\!\mu_i^2\! +\! 1\!\right)^2} + \frac{\nu_i^4 }{\left(\!\nu_i^2\! +\! 1\!\right)^2} \! \right)\! (F_i \!\circ\! F_i) 
& \left(\!\frac{\mu_i^2}{\left(\!\mu_i^2\! +\! 1\!\right)^2} + \frac{\nu_i^2}{\left(\nu_i^2\! +\! 1\right)^2} \!\right)\!(F_i\!\circ\! F_i)\\  
\left(\!\frac{\mu_i^2}{\left(\!\mu_i^2\! +\! 1\!\right)^2}+ \frac{\nu_i^2}{\left(\nu_i^2\! +\! 1\right)^2}\!\right)\!(F_i\!\circ \!F_i) 
&  \left(\!\frac{1}{\left(\!\mu_i^2\! +\! 1\!\right)^2} + \frac{1}{\left(\nu_i^2\! +\! 1\right)^2}\!\right) \!(F_i\! \circ \!F_i) } .
\]
Observe that $\mu_i + \nu_i = -\lambda_i$ and $\mu_i\nu_i = -1$, and so
\[
\mu_i^2 + \nu_i^2 = (\mu_i+\nu_i)^2- 2\mu_i\nu_i = \lambda_i^2 +2  
\]
and 
\[
\mu_i^4 + \nu_i^4 = \lambda_i^4 + 4\mu_i^2 + 4\nu_i^2 - 6 = \lambda_i^4 + 4 \lambda_i^2 + 2 .
\]We obtain that 
\[
\begin{split} 
\frac{\mu_i^4 }{\left(\mu_i^2\! +\! 1\right)^2} + \frac{\nu_i^4 }{\left(\nu_i^2\! +\! 1\right)^2}  &= \frac{\mu_i^4\left(\nu_i^2\! +\! 1\right)^2 + \nu_i^4  \left(\mu_i^2\! +\! 1\right)^2}{\left(\mu_i^2\! +\! 1\right)^2\left(\nu_i^2\! +\! 1\right)^2}\\
&= \frac{\lambda_i^4 + 6\mu_i^2 +6\nu_i^2 -4}{\left(\lambda_i^2 + 4\right)^2} 
= \frac{\lambda_i^2 + 2}{\lambda_i^2 + 4} 
\end{split}
\]
and 
\[
\frac{\mu_i^2 }{\left(\mu_i^2\! +\! 1\right)^2} + \frac{\nu_i^2 }{\left(\nu_i^2\! +\! 1\right)^2}  = \frac{\mu_i^2\left(\nu_i^2\! +\! 1\right)^2 + \nu_i^2  \left(\mu_i^2\! +\! 1\right)^2}{\left(\mu_i^2\! +\! 1\right)^2\left(\nu_i^2\! +\! 1\right)^2}
= \frac{2\lambda_i^2 +8}{\left(\lambda_i^2 + 4\right)^2} 
= \frac{2}{\lambda_i^2 + 4} 
\]
and 
\[
\frac{1 }{\left(\mu_i^2\! +\! 1\right)^2} + \frac{1 }{\left(\nu_i^2\! +\! 1\right)^2}  = \frac{\left(\nu_i^2\! +\! 1\right)^2 +  \left(\mu_i^2\! +\! 1\right)^2}{\left(\mu_i^2\! +\! 1\right)^2\left(\nu_i^2\! +\! 1\right)^2}
= \frac{\lambda_i^2 + 2}{\lambda_i^2 + 4} .
\]
Thus
\[
E_{\mu_i} \circ E_{\mu_i} + E_{\nu_i} \circ E_{\nu_i} =  
\pmat{\left(\frac{\lambda_i^2 + 2}{\lambda_i^2 + 4}  \right) (F_i \circ F_i) 
& \left(\frac{2}{\lambda_i^2 + 4}\right)(F_i\circ F_i)\\  
\left(\frac{2}{\lambda_i^2 + 4} \right)(F_i\circ F_i) 
&  \left(\frac{\lambda_i^2 + 2}{\lambda_i^2 + 4}\right) (F_i \circ F_i) } .
\]
Let \[ N =  \sum_{i=1}^n  \left(\frac{2}{\lambda_i^2 + 4}\right)(F_i\circ F_i).
\]
We can see that \[
\mxm(X(K_2)) = \pmat{\mxm(X) - N & N \\ N & \mxm(X) - N}
\]
and the lemma follows. \qed 

We observe that $N$ and $\mxm$ have the same kernel. 






\section{Trees with simple eigenvalues with $\rk(\mxm)$ bounded away from $\lceil \nicefrac{n}{2} \rceil$ } \label{sec:upperbd}

The following lemma is found in \cite{CoGdGuZh17}. 

\begin{lemma}\cite{CoGdGuZh17} If $X$ is a bipartite graph on $n$ vertices with $n$ distinct eigenvalues, then $rk(\mxm(X)) \leq \lceil \nicefrac{n}{2} \rceil$. \qed \end{lemma}

\begin{figure}[htbp]
    \centering
    \includegraphics{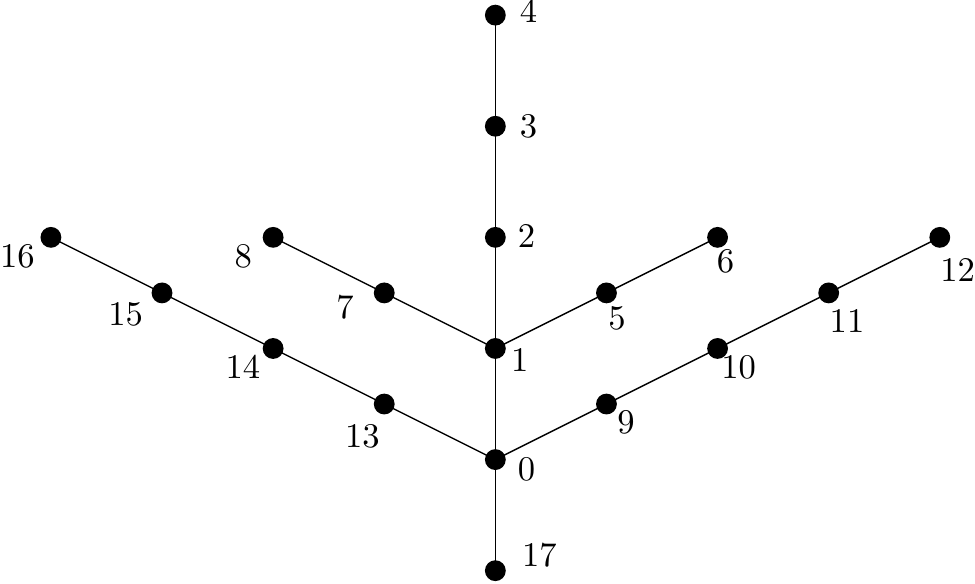}
    \caption{The tree $T^*$ on $18$ vertices such that $\mxm(T^*) = 8$. \label{fig:thetree}}
\end{figure}

In this section, we will give a construction for trees on $n$ vertices with simple eigenvalues where $\rk(\mxm)$ is bounded away from $\lceil \nicefrac{n}{2} \rceil$. A computer search finds that for $n = 1,\ldots, 17, 19,20$, every trees on $n$ vertices with simple eigenvalues has $\rk(\mxm) = \lceil \nicefrac{n}{2} \rceil$. Up to isomorphism, there is one tree $T^*$ on $18$ vertices such that $\rk(\mxm(T^*)) = 8$, which is given in Figure \ref{fig:thetree}. We will retain the notation of $T^*$ to denote this for the rest of this section.

We will need the following technical lemma. 

\begin{lemma} \label{lem:kerM}
Let $X$ be a graph with adjacency matrix  $A$. Suppose the spectral decomposition $A(X)= \sum_{r=0}^d \theta_r E_r$. Let $\mxm$ be the average mixing matrix of $X$. We have that 
\[  \ker(\mxm) = \bigcap_{r=0}^d \ker (E_r \circ E_r).\]
\end{lemma}

\proof Let $v \in \ker(\mxm)$. We see that 
\[
0= \Zv^TM\Zv = \sum_{r=0}^d \Zv^T E_r \circ E_r \Zv.
\]
Since $E_r \circ E_r$ is a principal submatrix of $E_r \otimes E_r$ and is hence positive semi-definite, we see that 
\[
\Zv^T E_r \circ E_r \Zv \geq 0
\]
for all $r$. Since the sum of non-negative real number is equal to zero, every term must be zero and thus
\[  \ker(\mxm) \subseteq \bigcap_{r=0}^d \ker (E_r \circ E_r).\]
Containment in the other direction is clear. \qed 

We will use Lemma \ref{lem:kerM} to bound the rank of iterated rooted products of $T^*$ with $K_2$. We note that Lemma \ref{lem:kerM} implies that vertices $u$ and $v$ are strongly cospectral in $X$ if and only if $\mxm e_u = \mxm e_v$, which is also proven in \cite{GodsilAverageMixing}. 

The characteristic polynomial of $T^*$ is 
\[
\begin{split}
\phi(T,x) = &(x - 1)  (x + 1) \left(x^2 - x - 1\right)\left(x^2 + x - 1]\right)\left(x^3 - x^2 - 2x + 1\right) \\
&\left(x^3 + x^2 - 2x - 1\right)\left(x^6 - 8x^4 + 12x^2 - 1\right)
\end{split}
\]
which has all simple roots. 

Let $X_0 = T^*$ and $X_{i+1} = X_i(K_2)$ for $i >0$. 

\begin{lemma} For $i\geq 0$, the graph $X_i$ has simple eigenvalues and
\[ \rk(\mxm(X_i)) \leq 2^{i+3}. \]\end{lemma}
 
 \proof We proceed by induction on $i$. We computed that $\rk(\mxm(T^*)) = 8$ and $T^*$ has simple eigenvalues. Suppose the statement is true for $i >1$. By Lemma \ref{lem:rootedsimp}, the graph $X_{i+1}$ has simple eigenvalues. For simplicity, let $n = 2^i18$, the number of vertices of $X_i$. Let $F_1,\ldots,F_n$ be  the orthogonal projections onto the eigenspaces of $A(X_i)$ with corresponding eigenvalues $\lambda_1,\ldots, \lambda_n$. Theorem \ref{thm:rootedmhat} gives that
\[
\mxm(X_{i+1}) = \pmat{\mxm(X_i) - N & N \\ N & \mxm(X_i) - N}
\]
where
\[ N =  \sum_{j=1}^n  \left(\frac{2}{\lambda_j^2 + 4}\right)(F_j\circ F_j).
\]
Let $\Zv \in \ker(\mxm(X_i))$. By Lemma \ref{lem:kerM}, we see that $\Zv \in \ker(F_j\circ F_j)$ for all $j$, and thus $N\Zv = 0$. Let $0_m$ denote the $m$ dimensional all zero vector.  We have 
\[
\mxm(X_{i+1}) \pmat{\Zv \\ 0_n} = \pmat{\mxm(X_i) - N & N \\ N & \mxm(X_i) - N}\pmat{\Zv \\ 0_n} = 0_{2n}
\]
and 
\[
\mxm(X_{i+1}) \pmat{0_n \\ \Zv} = \pmat{\mxm(X_i) - N & N \\ N & \mxm(X_i) - N}\pmat{0_n \\ \Zv}  = 0_{2n}.
\]
Thus, if $\{\Zv_k\}_{k}$ is an orthogonal basis for $\ker(\mxm(X_i))$, then 
\[
 \left\{ \pmat{\Zv_k \\ 0_n},\pmat{0_n \\ \Zv_k} \right\}_{k}
\]
is a set of $2 \dim \ker(\mxm(X_i))$ vectors in $\ker(\mxm(X_{i+1}))$ and the statement follows.\qed 

\begin{theorem}
For every positive real number $c$, there exists a tree $T$ with simple eigenvalues such that 
\[
\left\lceil \frac{|V(T)|}{2} \right\rceil -  \rk(\mxm(T)) > c.
\]
\end{theorem}

\proof Note that $X_i$ as construct above is a tree on $2^i18$ vertices such that 
\[
\left\lceil  \frac{|V(X_i)|}{2} \right\rceil  -  \rk(\mxm(X_i)) \geq 2^{i-1}18 - 2^i8 = 2^{i}. 
\]
For any $c>0$, we pick $i$ such that $2^i >c$ and the result follows. \qed 
 
\section{Lower bound on rank}\label{sec:lowerbd}
 
We denote by $X\diff u$ the graph obtaine from $X$ by deleting vertex $u$. The \textsl{coefficient matrix} of a graph $X$ is the matrix $C$ with rows indexed by the vertices of $X$ and columns indexed by integers $[n]$ such  that the $(u,r)$ entry of $C$ is the coefficient of $t^{r-1}$ in the characteristic polynomial $\phi(X\diff u,t)$. We will see that for graphs with simple eigenvalues, the coefficient matrix is useful in studying the average mixing matrix. The following lemma appears in the proof of \cite[Theorem 3.2]{GodsilAverageMixing}. 

\begin{lemma}\cite{GodsilAverageMixing} Let $X$ be a graph on $n$ vertices with $n$ distinct eigenvalues $\theta_1, \ldots, \theta_n$. Let $\Delta$ be the $n × n$ diagonal matrix whose $r$-th diagonal entry is $\phi'(X,\theta_r)$. Let $V$ be the
$n ×n$ Vandermonde matrix with $ij$-entry equal to $\theta^{i-1}_j$. Then 
\[
\mxm(X) = CV\Delta^{-2}V^TC^T. 
\] \qed 
\end{lemma}
 
 From this, we derive the following corollary as a direct consequence. 
 
\begin{corollary}\label{lem:coef}
If $X$ is a graph with simple eigenvalues, then $\rk ( \mxm(X))$ is equal to the rank of the coefficient matrix. 
\end{corollary}

We will also use this standard fact, which can be found in \cite{GodsilAlgebraicCombinatorics}, for example. 
 
\begin{lemma}\label{lem:charpoly=matchpoly}
If $T$ is a tree, then the characteristic polynomial of $T$ is equal to the matching polynomial of $T$. \qed
\end{lemma}

Lemma~\ref{lem:submatrixfullrank} is Corollary 8.9.2 in \cite{GR}. 

\begin{lemma}\cite{GR}\label{lem:submatrixfullrank} If $A$ is a symmetric matrix of rank $r$, then it has a principal submatrix of full rank. \qed
\end{lemma}

Corollary~\ref{cor:perfect} follows from Lemma~\ref{lem:submatrixfullrank} and the fact that a tree has a perfect matching if and only if its adjacency matrix has full rank.

\begin{corollary}\label{cor:perfect} If $T$ is a tree with simple eigenvalues, then either $T$ has a perfect matching or there exists a vertex $v$ such that $T\diff v$ has a perfect matching.
\end{corollary}

\proof 
If $T$ is a tree with simple eigenvalues, then zero can have multiplicity at most one.  If the multiplicity of zero is zero, then $T$ has a perfect matching.  If the multiplicity of zero is one, then the rank of $A(T)$ is $|T|-1$.  Thus, by Lemma~\ref{lem:submatrixfullrank}, $A(T)$ has $(|T|-1)\times (|T|-1)$ principal submatrix of full rank.  Thus, there exists a vertex $v$ such that $A(T\diff v)$ has full rank.  It follows that $T\diff v$ has a perfect matching.  
\qed 

\begin{lemma}\label{lem:degtwo}
If $T$ is a tree with simple eigenvalues, $|T|\geq 3$, then $T$ has a leaf adjacent to a vertex of degree two.
\end{lemma}

\proof 
 The path on three vertices satisfies the conclusion of the lemma.  Let $T$ be a tree with simple eigenvalues and $|T|\geq 4$. Consider a diametrical path $P=\{v_0,\ldots, v_d\}$ of $T$.  Since the star on $n\geq 4$ vertices does not have simple eigenvalues, $d\geq 3$.  

Since $P$ is a diametrical path of length at least three, $v_1$ and $v_{d-1}$ are distinct and adjacent to exactly one non leaf vertex of $T$.  If either $v_1$ or $v_{d-1}$ has degree two, then we are done.  Suppose by way of contradiction, that both vertices have degree greater than two.  Deleting $v_1$ and $v_{d-1}$ from $T$ yields a graph with at least four isolated vertices. As each isolated vertex contributes a zero eigenvalue, the multiplicity of zero is at least four.  Since the eigenvalues of induced subgraphs interlace, the multiplicity of the zero eigenvalue of $T$ is at least 2. This contradicts that $T$ has simple eigenvalues. \qed

\begin{theorem}
If $T$ is a tree with simple eigenvalues on at least four vertices and is not isomorphic to $P_4$, then $\rk\left(\mxm(T)\right)\geq 3$.
\end{theorem}
 
\proof In light of Corollary~\ref{lem:coef}, it is sufficient to show that the rank of the coefficient matrix of a tree is at least three.  By Lemma~\ref{lem:degtwo}, $T$ has a leaf $u$ adjacent to a vertex $v$ of degree two.  Let $w$ be the non leaf neighbor of $v$ with $\deg(w)=\ell$.  Let $m_{\alpha}(T)$ denote the number of matchings of size $\alpha$ in $T$. By Lemma~\ref{lem:charpoly=matchpoly}, the $i$th row of the coefficient matrix of $T$ can be expressed in terms of $m_\alpha(T\diff i)$ for the appropriate value of $\alpha$.
\medskip

\noindent \textbf{Case 1:}  $T$ has a perfect matching.

\medskip
As $T$ has a perfect matching, $n=|T|$ is even. The only tree on four vertices with a perfect matching is $P_4$.  It has two pairs of cospectral vertices and thus $\rk (\mxm(P_4))=2$. For the remainder of this case, assume that $n\geq 6$.  Consider the $3\times 3$ submatrix $C_1$ of the coefficient matrix of $T$ with rows corresponding to $u,v,w$ and columns corresponding to $t, t^{n-3}$, and $t^{n-1}$. For simplicity in the notation let $k=\nicefrac{n}{2}$.  Then

\[C_1= 
\pmat{
(-1)^{k-1}m_{k-1}(T\diff u) & -m_1(T\diff u) & m_0(T\diff u)\\
(-1)^{k-1}m_{k-1}(T\diff v) & -m_1(T\diff v) & m_0(T\diff v)\\
(-1)^{k-1}m_{k-1}(T\diff w) & -m_1(T\diff w) & m_0(T\diff w)\\
}
\]

As $m_0(T\diff i)$ is defined to be 1 and $m_1(T\diff i)$ is the number of edges in $T\diff i$, it remains to determine the values of the first column.  

In $T\diff u$, the set of matchings of size $k-1$ can be partitioned into two sets: the $k-1$ matchings which use edge $vw$ and those that do not.  Thus $$m_{k-1}(T\diff u)=m_{k-2}(T\diff \{u,v,w\})+m_{k-1}(T\diff \{u,v\}).$$

Since $T$ has a perfect matching and every perfect matching of $T$ contains the edge $uv$, we have that $T\diff \{u,v\}$ also contains a perfect matching.  Thus $m_{k-1}(T\diff \{u,v\})=1$.  Similarly, $m_{k-1}(T\diff v)=1$.

Since $T\diff w$ has $uv$ as a component and any matching of size $k-1$ uses $uv$, we see that \[m_{k-1}(T\diff w)=m_{k-2}(T\diff \{u,v,w\}.\]  Letting $q=m_{k-2}(T\diff \{u,v,w\})$,  we obtain
\[C_1=\pmat{
(-1)^{k-1}(q+1) & -(n-2) & 1\\
(-1)^{k-1} & -(n-3) & 1\\
(-1)^{k-1}q & -(n-\ell-1) & 1\\
} \]
and $$\det C_1=(-1)^{k-1}(1+q(1-\ell)).$$
Note that $\det C_1=0$ if and only if $\ell=1+\nicefrac{1}{q}$.  As $\ell\in\mathbb{N}$, this occurs if and only if $q=1$ and  $\ell=2$. We claim that if $\ell=2$, then $q\geq 2$. Let $y$ be the other neighbor of $w$.  Since $T$ has a perfect matching, so does $T\diff \{u,v,w,y\}$.  Thus $m_{k-2}(T\diff \{u,v,w,y\})=1$ and $T\diff \{u,v,w\}$ has a $k-2$ matching which uses no edge incident to $y$.  Since $n\geq 6$, $\deg(y)\geq 2$.  Thus $T\diff \{u,v,w\}$ has at least one $k-2$ matching which uses an edge incident to $y$.  Therefore $q=m_{k-2}(T\diff \{u,v,w\}\geq 2$.  Further $\det C_1\neq 0$ and the coefficient matrix of $T$ has rank at least 3.
\medskip

\noindent \textbf{Case 2:}  $T$ does not have a perfect matching.

\medskip
There are two possibilities to consider;  either $T\diff u$ has a perfect matching or it does not.  

In the case that $T\diff u$ has a perfect matching, $n=|T|$ is odd and so $|T|\geq 5$. Consider the $3\times 3$ submatrix $C_2$ of the coefficient matrix of $T$ with rows corresponding to $u,v,w$ and columns corresponding to $t^0, t^{n-3}$, and $t^{n-1}$. For simplicity in the notation let $\nicefrac{(n-1)}{ 2}=j$.  Then

\[C_2= 
\pmat{
(-1)^{j}m_{j}(T\diff u) & -m_1(T\diff u) & m_0(T\diff u)\\
(-1)^{j}m_{j}(T\diff v) & -m_1(T\diff v) & m_0(T\diff v)\\
(-1)^{j}m_{j}(T\diff w) & -m_1(T\diff w) & m_0(T\diff w)\\
}.\]

As stated in Case 1, $m_0(T\diff i)=1$ and $m_1(T\diff i)$ is  the number of edges in $T\diff i$.  Matchings of size $j=\nicefrac{n-1}{2}$ are perfect matchings in a graph on $n-1$ vertices.  Further since $T\diff u$ has a perfect matching so does $T\diff w$.  Lastly, $T\diff v$ does not have a perfect matching as $u$ is isolated by deleting $v$.  Thus 
\[C_2=\pmat{
(-1)^{j} & -(n-2) & 1\\
0 & -(n-3) & 1\\
(-1)^{j} & -(n-\ell-1) & 1\\
}.\]
Note that $\det C_2=(-1)^j(\ell-1)$.  Since $\deg(w)=\ell\geq 2$, $C_2$ has full rank and so the coefficient matrix of $T$ has rank at least three.

We now consider the other possibility mentioned at the start of this case.  Assume that $T\diff u$ does not have a perfect matching.  By Corollary~\ref{cor:perfect} there exists a vertex $z$ of degree $f$ such that $T-z$ has a perfect matching.  Thus $n=|T|$ is odd and $n\geq 5$.  Consider the $3\times 3$ submatrix $C_3$ of the coefficient matrix of $T$ with rows corresponding to $u,v,z$ and columns corresponding to $t^0, t^{n-3}$, and $t^{n-1}$. For simplicity in the notation let $\nicefrac{n-1}{2}=j$.  Then
\[C_3= 
\pmat{
(-1)^{j}m_{j}(T\diff u) & -m_1(T\diff u) & m_0(T\diff u)\\
(-1)^{j}m_{j}(T\diff v) & -m_1(T\diff v) & m_0(T\diff v)\\
(-1)^{j}m_{j}(T\diff z) & -m_1(T\diff z) & m_0(T\diff z)\\
}.\]

As in the other cases, the entries of the last two columns are easily filled. The entries in the first column are determined by whether $T\diff i$ has a perfect matching.  Thus,
\[C_3= 
\pmat{
0 & -(n-2) & 1\\
0 & -(n-3) & 1\\
(-1)^{j} & -(n-f-1) & 1\\
}. \]
Note that $\det C_3=(-1)^{j+1}\neq 0$. Thus $C_3$ has full rank and so the coefficient matrix of $T$ has rank at least three. \qed 

\section{Computational data} 

Computations of average mixing matrices of trees on up to $20$ were carried out, using Sage Mathematics Software \cite{sage}. The results are recorded in Tables \ref{tab:trees2-14} and \ref{tab:trees15-20}. For each order $n$ and rank $r$, we record the number of trees on $n$ vertices whose average mixing matrix has rank $r$ and, amongst those, the number with $n$ distinct eigenvalues.

\begin{table}[htbp]
\centering
\begin{tabular}{|llll|}
\hline
n  & rank & \# trees & \# simple  \\
&        & & eigenvalues \\
\hline
2  & 1    & 1        & 1                     \\
\hline
3  & 2    & 1        & 1                     \\
\hline
4  & 2    & 1        & 1                     \\
   & 4    & 1        & 0                     \\
\hline
5  & 3    & 2        & 2                     \\
   & 5    & 1        & 0                     \\
\hline
6  & 3    & 3        & 2                     \\
   & 5    & 2        & 0                     \\
   & 6    & 1        & 0                     \\
\hline
7  & 4    & 5        & 5                     \\
   & 5    & 1        & 0                     \\
   & 6    & 4        & 0                     \\
   & 7    & 1        & 0                     \\
\hline
8  & 4    & 5        & 4                     \\
   & 5    & 4        & 0                     \\
   & 6    & 8        & 0                     \\
   & 7    & 4        & 0                     \\
   & 8    & 2        & 0                     \\
\hline
9  & 5    & 19       & 18                    \\
   & 6    & 3        & 0                     \\
   & 7    & 15       & 0                     \\
   & 8    & 7        & 0                     \\
   & 9    & 3        & 0                     \\
\hline
10 & 4    & 1        & 0                     \\
   & 5    & 14       & 11                    \\
   & 6    & 19       & 0                     \\
   & 7    & 30       & 0                     \\
   & 8    & 21       & 0                     \\
   & 9    & 16       & 0                     \\
   & 10   & 5        & 0                     \\
\hline
\end{tabular}
\hspace{15pt}
\begin{tabular}{|llll|}
\hline
n  & rank & \# trees & \# simple  \\
&        & & eigenvalues \\
\hline
11 & 5    & 1        & 0                     \\
   & 6    & 64       & 62                    \\
   & 7    & 18       & 0                     \\
   & 8    & 79       & 0                     \\
   & 9    & 40       & 0                     \\
   & 10   & 26       & 0                     \\
   & 11   & 7        & 0                     \\
\hline
12 & 5    & 1        & 0                     \\
   & 6    & 44       & 37                    \\
   & 7    & 106      & 0                     \\
   & 8    & 129      & 0                     \\
   & 9    & 119      & 0                     \\
   & 10   & 93       & 0                     \\
   & 11   & 48       & 0                     \\
   & 12   & 11       & 0                     \\
\hline
13 & 6    & 2        & 0                     \\
   & 7    & 264      & 250                   \\
   & 8    & 107      & 0                     \\
   & 9    & 411      & 0                     \\
   & 10   & 223      & 0                     \\
   & 11   & 186      & 0                     \\
   & 12   & 87       & 0                     \\
   & 13   & 21       & 0                     \\
\hline
14 & 6    & 4        & 0                     \\
   & 7    & 146      & 116                   \\
   & 8    & 552      & 0                     \\
   & 9    & 591      & 0                     \\
   & 10   & 694      & 0                     \\
   & 11   & 622      & 0                     \\
   & 12   & 341      & 0                     \\
   & 13   & 172      & 0                     \\
   & 14   & 37       & 0                     \\
\hline
\end{tabular}
\caption{Ranks of average mixing matrices of trees on $2$ to $14$ vertices. \label{tab:trees2-14}}
\end{table}

\begin{table}[htbp]
\centering
\begin{tabular}{|llll|}
\hline
n  & rank & \# trees & \# simple  \\
&        & & eigenvalues \\
\hline
15 & 7    & 4        & 0                     \\
   & 8    & 1117     & 1041                  \\
   & 9    & 663      & 0                     \\
   & 10   & 2173     & 0                     \\
   & 11   & 1365     & 0                     \\
   & 12   & 1328     & 0                     \\
   & 13   & 719      & 0                     \\
   & 14   & 309      & 0                     \\
   & 15   & 63       & 0                     \\
\hline
16 & 7    & 7        & 0                     \\
   & 8    & 543      & 465                   \\
   & 9    & 2926     & 0                     \\
   & 10   & 2834     & 0                     \\
   & 11   & 4265     & 0                     \\
   & 12   & 3881     & 0                     \\
   & 13   & 2650     & 0                     \\
   & 14   & 1494     & 0                     \\
   & 15   & 600      & 0                     \\
   & 16   & 120      & 0                     \\
\hline
17 & 8    & 11       & 0                     \\
   & 9    & 4889     & 4452                  \\
   & 10   & 4325     & 0                     \\
   & 11   & 11653    & 0                     \\
   & 12   & 8340     & 0                     \\
   & 13   & 9347     & 0                     \\
   & 14   & 5724     & 0                     \\
   & 15   & 3002     & 0                     \\
   & 16   & 1146     & 0                     \\
   & 17   & 192      & 0                     \\
\hline
18 & 7    & 2        & 0                     \\
   & 8    & 25       & 1                     \\
   & 9    & 2108     & 1727                  \\
   & 10   & 15306    & 0                     \\
\hline
\end{tabular}
\hspace{15pt}
\begin{tabular}{|llll|}
\hline
n  & rank & \# trees & \# simple  \\
&        & & eigenvalues \\
\hline
18   & 11   & 14829    & 0                     \\
   & 12   & 26545    & 0                     \\
   & 13   & 24194    & 0                     \\
   & 14   & 19249    & 0                     \\
   & 15   & 12980    & 0                     \\
   & 16   & 6019     & 0                     \\
   & 17   & 2242     & 0                     \\
   & 18   & 368      & 0                     \\
\hline
19 & 8    & 2        & 0                     \\
   & 9    & 25       & 0                     \\
   & 10   & 22159    & 19884                 \\
   & 11   & 26204    & 0                     \\
   & 12   & 64701    & 0                     \\
   & 13   & 53492    & 0                     \\
   & 14   & 63220    & 0                     \\
   & 15   & 43183    & 0                     \\
   & 16   & 27389    & 0                     \\
   & 17   & 12603    & 0                     \\
   & 18   & 4259     & 0                     \\
   & 19   & 718      & 0                     \\
\hline
20 & 8    & 5        & 0                     \\
   & 9    & 43       & 0                     \\
   & 10   & 8641     & 7055                  \\
   & 11   & 81498    & 0                     \\
   & 12   & 79080    & 0                     \\
   & 13   & 165082   & 0                     \\
   & 14   & 153019   & 0                     \\
   & 15   & 139556   & 0                     \\
   & 16   & 102182   & 0                     \\
   & 17   & 58113    & 0                     \\
   & 18   & 26098    & 0                     \\
   & 19   & 8405     & 0                     \\
   & 20   & 1343     & 0                    \\
\hline
\end{tabular}
\caption{Ranks of average mixing matrices of trees on $15$ to $20$ vertices. \label{tab:trees15-20}}
\end{table}

\section{Open problems}

In Section \ref{sec:lowerbd}, we give a constant lower on the rank of the average mixing matrix of trees with simple eigenvalues on at least $4$ vertices. Table \ref{tab:minrks} shows the minimum ranks among the average mixing matrices of trees of $n$ vertices. 

\begin{table}[H]
\centering
\renewcommand\arraystretch{1.1}
\renewcommand\tabcolsep{4pt}
\begin{tabular}{|c|lllllllllllllllllll|}
\hline
$n$        & 2 & 3 & 4 & 5 & 6 & 7 & 8 & 9 & 10 & 11 & 12 & 13 & 14 & 15 & 16 & 17 & 18 & 19 & 20 \\
\hline 
min rank & 1 & 2 & 2 & 3 & 3 & 4 & 4 & 5 & 4  & 5  & 5  & 6  & 6  & 7  & 7  & 8  & 7  & 8  & 8 \\
\hline
\end{tabular}
\caption{The minimum ranks of the average mixing matrices of trees on $n$ vertices, for $n = 2, \ldots, 20$.  \label{tab:minrks}}
\end{table}

We are motivated by the computational data to ask the following question:

\begin{question} Does there exist a non-constant, increasing function $f(n)$ such that 
\[\rk(\mxm(T)) \geq f(n) \]
for any tree $T$ on $n$ vertices?
\end{question}


\begin{thebibliography}{10}

\bibitem{TamonAdamczakUniformMixingCycles}
William Adamczak, Kevin Andrew, Leon Bergen, Dillon Ethier, Peter Hernberg,
  Jennifer Lin, and Christino Tamon.
\newblock {Non-Uniform Mixing Of Quantum Walk On Cycles}.
\newblock {\em International Journal Of Quantum Information}, 05(06):781--793,
  2007.

\bibitem{TamonAdamzackAverageMixing}
William Adamczak, Kevin Andrew, Peter Hernberg, and Christino Tamon.
\newblock {A note on graphs resistant to quantum uniform mixing}.
\newblock {\em ArXiv e-prints}, page~9, 08 2003.

\bibitem{AdaChanComplexHadamardIUMPST}
Ada Chan.
\newblock {Complex Hadamard Matrices, Instantaneous Uniform Mixing and Cubes}.
\newblock {\em ArXiv e-prints}, 2013.

\bibitem{ChildsUniversalQComputation}
Andrew~M Childs.
\newblock {Universal computation by quantum walk}.
\newblock {\em Physical Review Letters}, 102(18):4,180501, 2009.

\bibitem{CoutinhoGodsilGuoVanhove2}
G.~Coutinho, C.~Godsil, K.~Guo, and F.~Vanhove.
\newblock Perfect state transfer on distance-regular graphs and association
  schemes.
\newblock {\em Linear Algebra and its Applications}, 478:108 -- 130, 2015.

\bibitem{CoGdGuZh17}
G.~{Coutinho}, C.~{Godsil}, K.~{Guo}, and H.~{Zhan}.
\newblock {A New Perspective on the Average Mixing Matrix}.
\newblock {\em ArXiv e-prints}, September 2017.

\bibitem{FarhiGutmann}
Edward Farhi and Sam Gutmann.
\newblock {Quantum computation and decision trees}.
\newblock {\em Physical Review A}, 58(2):915--928, 1998.

\bibitem{GR}
C.~Godsil and G.~Royle.
\newblock {\em Algebraic Graph Theory}, volume 207 of {\em Graduate Texts in
  Mathematics}.
\newblock Springer-Verlag, New York, 2001.

\bibitem{GoMc78}
C.D. Godsil and B.D. McKay.
\newblock A new graph product and its spectrum.
\newblock {\em Bulletin of the Australian Mathematical Society}, 18(1):21–28,
  1978.

\bibitem{GodsilAlgebraicCombinatorics}
Chris~D Godsil.
\newblock {\em {Algebraic Combinatorics}}.
\newblock Chapman {\&} Hall, New York, 1993.

\bibitem{GodsilStateTransfer12}
Chris~D Godsil.
\newblock {State transfer on graphs}.
\newblock {\em Discrete Mathematics}, 312(1):129--147, 2012.

\bibitem{GodsilAverageMixing}
Chris~D Godsil.
\newblock {Average mixing of continuous quantum walks}.
\newblock {\em Journal of Combinatorial Theory, Series A}, 120(7):1649--1662,
  2013.

\bibitem{GodsilMullinRoy}
Chris~D Godsil, Natalie Mullin, and Aidan Roy.
\newblock {Uniform Mixing and Association Schemes}.
\newblock {\em ArXiv e-prints}, 2013.

\bibitem{KayReviewPST}
Alastair Kay.
\newblock {Perfect, efficient, state transfer and its application as a
  constructive tool}.
\newblock {\em International Journal of Quantum Information}, 8(04):641--676,
  2010.

\bibitem{sage}
W.\thinspace{}A. Stein et~al.
\newblock {\em {S}age {M}athematics {S}oftware ({V}ersion 6.1.1)}.
\newblock The Sage Development Team, 2014.
\newblock {\tt http://www.sagemath.org}.

\bibitem{VinetZhedanovHowTo}
Luc Vinet and Alexei Zhedanov.
\newblock {How to construct spin chains with perfect state transfer}.
\newblock {\em Physical Review A}, 85(1):12323, 2012.

\end{thebibliography}

\end{document}